\documentclass[12pt]{amsart}
\usepackage{amssymb}

\newtheorem{Thm}{Theorem}[section]

\newtheorem{Cor}[Thm]{Corollary}
\newtheorem{Lem}[Thm]{Lemma}
\newtheorem{Prop}[Thm]{Proposition}

\def\ldots{\mathinner{\ldotp\ldotp\ldotp}}
\def\ldots{\mathinner{\cdotp\cdotp\cdotp}}

\begin{document}

\title{ The $H^{\infty}-$calculus and sums of closed operators
}

\author{N. J. Kalton}
\address{Department of Mathematics \\
University of Missouri-Columbia \\
Columbia, MO 65211\\ USA
}

\email{nigel@math.missouri.edu}
\author{L. Weis}\address{Mathematisches Institut I\\ Universit\"at Karlsruhe\\ 76128 Karlsruhe\\ Germany}  \email{Lutz.Weis@math.uni-Karlsruhe.de}

\subjclass{Primary: 47A60, 47D06}
\thanks{The first author acknowledges support from NSF grant DMS-9800027
  and from DFG grant Ka 97/126.}

\begin{abstract} We develop a very general operator-valued functional calculus
  for operators with an $H^{\infty}-$calculus. We then apply this to the joint
  functional calculus of two commuting sectorial operators when one has an
  $H^{\infty}$calculus.  Using this we prove theorem of Dore-Venni type on
  sums of commuting sectorial operators and apply our results to the problem
  of $L_p-$maximal regularity.
  Our main assumption is the R-boundedness of certain sets of operators, and
  therefore methods from the geometry of Banach spaces are essential here.
  In the final section we exploit the special Banach space structure of
  $L_1-$spaces and $C(K)-$spaces, to obtain some more detailed results in this
  setting.
\end{abstract}

\maketitle

\section{Introduction}\label{intro}

In recent years the  notion of an $H^{\infty}-$calculus  for
sectorial  operators  on  a  Banach  space  has  played   an
important role  in spectral  theory for  unbounded operators
and its applications to differential operators and evolution
equations.    We  recall  that  a sectorial operator of type
$0\le \omega<\pi$ satisfies a ``parabolic'' estimate of  the
type \begin{equation}\label{parabolic} \|\zeta R(\zeta, A)\|
\le C_{\sigma}\qquad  |\arg \zeta|\ge  \sigma \end{equation}
for  every  $\omega<\sigma<\pi.$  This  estimate  allows   a
definition of $f(A)$ as a bounded operator for functions $f$ which are bounded and
analytic    on    the    sector   $\Sigma_{\sigma}=\{\zeta:\
|\arg\zeta|<\sigma\}$ and which obey a condition of the type
$|f(\zeta)|\le C(|\zeta|/(1+|\zeta|^2)^{\epsilon})$ for  some
$\epsilon>0.$ This  is described  in detail  in \cite{M} and
\cite{LLLeM} and  we give  a somewhat  different approach in
Section \ref{functional} below.  If we then have an estimate
$$\|f(A)\| \le  C\|f\|_{H^{\infty}(\Sigma_{\sigma})}$$ it  is
possible to  extend the  definition of  $f(A)$ to  any $f\in
H^{\infty}(\Sigma_{\sigma})$  and  we  say  that  $f$ has an
$H^{\infty}(\Sigma_{\sigma})-$calculus.    It  is,  by  now,
well-known  that  many  systems  of  parabolic  differential
operators, Schr\"odinger  operators and  pseudo-differential
operators do have an $H^{\infty}-$calculus (cf.  \cite{HP2},
\cite{AMN}, \cite{AHS}  and \cite{DM})  and this  has proved
useful in applications.

Of particular importance  are two closely  related problems:
\begin{itemize}\item  the  maximal  $L_p-$regularity  of the
Cauchy problem  $$ y'(t)+Ay(t)  =f(t),\qquad y(0)=0$$  for a
sectorial operator  of type  $\omega<\frac{\pi}2$ \item  the
question  whether  the  sum  $A+B$  with  domain   $\mathcal
D(A)\cap\mathcal D(B)$ of two commuting sectorial  operators
is closed.  \end{itemize} In  fact the first problem can  be
reduced to the second, and the latter problem is essentially
the question  whether one  can construct  a bounded operator
$B(A+B)^{-1}.$ This then is a special case of the problem of
constructing a  joint functional  calculus of  $A,B.$ In the
case of Hilbert  spaces and some  related situations it  was
shown in  \cite{FM}, \cite{LLLeM}  and \cite{LLeM}  that one
can   construct   an   operator-valued  functional  calculus
associated  to  an  operator  with $H^{\infty}-$calculus and
this permits  a solution;  however, it  was also  shown that
such an approach  cannot work in  general Banach spaces  and
additional conditions are therefore needed.

We now describe the main results of this paper.  In  Section
\ref{functional}  we  describe  a  method  of setting up the
joint  functional  calculus   of  $n$  commuting   sectorial
operators  and  an  operator-valued  extension.   In Section
\ref{Rademacher}     we     recall     the     notion     of
Rademacher-boundedness  (or  R-boundedness)  of  families of
operators.  This  implicitly goes back  to work of  Bourgain
\cite {Bo} and has recently been studied in \cite {BG}, \cite
{CPSW},   \cite{CP}   and   \cite{We}   in  connection  with
vector-valued multiplier theorems.   We also  introduce some
weaker  notions  and  study  their  relationship  to certain
Banach space properties of the underlying space.

Using these ideas in Section \ref{operator} we prove a  very
general  result  on  the  existence  of  an  operator-valued
functional     calculus     for     operators     with    an
$H^{\infty}$-calculus.  Roughly speaking this permits us  to
replace  boundedness  of  the  range  of  the  function   by
Rademacher-boundedness  (or  even  the  weaker  concept   of
U-boundedness introduced in Section \ref{Rademacher}).

We then study the relationship between R-boundedness and the
functional  calculus  for  general  sectorial operators.  Of
particular importance is  the notion of  R-sectoriality when
the boundedness condition  (\ref{parabolic}) is replaced  by
an  R-boundedness   condition.     Using  this   in  Theorem
\ref{sumthm} we prove a general result on sums which can  be
regarded  as   an  extension   of  the   Dore-Venni  Theorem
\cite{DV}.  We  show that if  $A,B$ are commuting  sectorial
operators       such       that       $A$       has       an
$H^{\infty}(\Sigma_{\sigma})-$calculus     and     $B$    is
R-sectorial of type $\sigma'$ when $\sigma+\sigma'<\pi$ then
$A+B$  (with  domain  $\mathcal  D(A)\cap\mathcal  D(B))$ is
closed.  One practical advantage of our result is that it is
easier  to  check  R-sectoriality  than  the boundedness of
imaginary powers (see \cite{Ku}, \cite{We2} and  \cite{CP}).
We     also     give     applications     to    the    joint
$H^{\infty}-$functional calculus (cf.  \cite{LLeM}) and show
how Banach space  properties such as  UMD, analytic UMD  and
property $(\alpha)$ of Pisier relate to our results.

It   might    be   added    that   our    results   on   the
$H^{\infty}-$calculus   emphasize    the   fact    that   an
$H^{\infty}-$calculus   really   induces   an  unconditional
expansion of  the identity  on the  underlying Banach space.
We  feel  our  development  of  the  theory here is somewhat
simpler than preceding approachs (even for Hilbert spaces).

Finally in  Section \ref{Groth}  we use  this observation to
show how  classical results  on unconditional  bases due  to
Lindenstrauss and Pe\l czy\'nski \cite{LP} can be recast  as
results on operators with an $H^{\infty}-$calculus on  $L_1$
and  $C(K)-$spaces.    In  these  cases  we  get very strong
conclusions, but they are mitigated by the observation  that
there are in practice very few examples of such operators on
spaces of this type.

\section{An   operator-valued   functional   calculus    for
sectorial                       operators}\label{functional}
\setcounter{equation}{0} In this section we sketch a  method
of setting up an  operator-valued functional calculus for  a
sectorial operator  and a  joint operator-valued  functional
calculus  for  finite  collections  of  commuting  sectorial
operators.     For  an   alternative  construction   of  the
$H^{\infty}$-calculus  based  on  McIntosh's  approach   see
\cite{LLeM}.

Let   us   first   introduce   some   notation.      Suppose
$0<\sigma<\pi.$  Then  we  denote  by  $\Sigma_{\sigma}$ the
sector   $\{z:      |\arg   z|<\sigma,   |z|>0\}$   and   by
$\Gamma_{\sigma}$     the     contour     $\{|t|e^{i(\text{sgn
}t)\sigma}:\    -\infty<t<\infty\}.$     We    denote     by
$H^{\infty}(\Sigma_{\sigma})$  the  space  of  all   bounded
analytic   functions   on   $\Sigma_{\sigma}.$   We   define
$H^{\infty}_0(\Sigma_{\sigma})$ to  be the  subspace of  all
$f\in H^{\infty}(\Sigma_{ \sigma}) $ which obey an  estimate
of  the  form  $|f(z)|\le  C(|z|/(1+|z|^2))^{\epsilon}$ with
$\epsilon>0.$  Let  us  extend  this  to  dimension $m$.  In
$\mathbb C^m$  if $\sigma=(\sigma_1,\ldots,\sigma_m)$  where
$0<\sigma_k<\pi$                  we                  define
$\Sigma_{\sigma}=\prod_{k=1}^m\Sigma_{\sigma_k}$         and
$\Gamma_{\sigma}=\prod_{k=1}^m\Gamma_{\sigma_k}$.         If
$\sigma,\nu\in  \mathbb  R^m$   we  write  $\sigma>\nu$   if
$\sigma_k>\nu_k$   for   $1\le   k\le   m.$   We  denote  by
$H^{\infty}(\Sigma_{\sigma})$  the  space  of  all   bounded
analytic   functions   on   $\Sigma_{\sigma}.$   We   define
$H^{\infty}_0(\Sigma_{\sigma})$ to  be the  subspace of  all
$f\in H^{\infty}(\Sigma_{ \sigma}) $ which obey an  estimate
of              the              form             $|f(z)|\le
C\prod_{k=1}^m(|z_k|/(1+|z_k|^2))^{\epsilon}$           with
$\epsilon>0$ where $z=(z_1,\ldots, z_m).$

Next we introduce  some corresponding vector-valued  spaces.
Now suppose $X$  is a Banach  space and ${\mathcal  A}$ is a
subalgebra  of  $\mathcal  L(X),$  which  is  closed for the
strong-operator           topology.                       If
$\sigma=(\sigma_1,\ldots,\sigma_m)$  as  above,  we   define
$H^{\infty}(\Sigma_{\sigma};\mathcal  A)$  the  space of all
bounded  functions  $F:\Sigma_{\sigma}\to{\mathcal  A},$  so
that for  every $x\in  X$ the  map $z\to  F(z)x$ is analytic
(i.e.   $F$ is  analytic for  the strong-operator topology).
We consider  the scalar  space $H^{\infty}(\Sigma_{\sigma})$
as a subspace of $H^{\infty}(\Sigma_{\sigma};{\mathcal  A})$
via the identification $f \to  fI.$ We shall say that  $F_n$
{\it      converges      boundedly}      to      $F$      in
$H^{\infty}(\Sigma_{\sigma};\mathcal         A)$          if
$\sup_n\sup_{z\in   \Sigma_{\sigma}}\|F_n(z)\|<\infty$   and
$F_n(z)x\to  F(z)x$  for  every  $z\in \Sigma_{\sigma},$ and
$x\in X.$  We define  $H^{\infty}_0(\Sigma_{\sigma},\mathcal
A)$   the   subspace   of   all   $F\in  H^{\infty}(\Sigma_{
\sigma},{\mathcal A}) $ which  obey an estimate of  the form
$\|F(z)\|\le    C\prod_{k=1}^m(|z_k|/(1+|z_k|^2))^{\epsilon}$
with $\epsilon>0$ where $z=(z_1,\ldots, z_m).$

We next consider the space of germs of such functions.   Fix
$0\le \omega_k<\pi$ for $1\le k\le m.$ We consider the space
$\mathcal          H(\omega,\mathcal          A)           =
\cup_{\sigma>\omega}H^{\infty}(\Sigma_{\sigma};{\mathcal
A}))$  where   $(F,G)$  are   identified  if   there  exists
$\sigma>\omega$    with    $F(z)=G(z)$    for    all   $z\in
\Sigma_{\sigma}.$ $\mathcal H(\omega,\mathcal A)$ is then an
algebra.   In $\mathcal  H(\omega,\mathcal A)$  we define  a
notion of sequential  convergence $\tau$ by  $F_n\to F $  if
there  exists   $\sigma>\omega$  so   that  each   $F_n,F\in
H^{\infty}(\Sigma_{\sigma};\mathcal  A)$,  $\sup_n\sup_{z\in
\Sigma_{\sigma}}\|F_n(z)\|<\infty$ and $F_n(z)x\to F(z)x$  for
all $z\in \Sigma_{\sigma}$
and all $x\in X.$

Recall  that  a  closed  densely  defined  operator $A$ on a
Banach  space  $X$  is a {\it  sectorial  operator  of   type
$0\le\omega=\omega(A)<\pi$}  if  $A$  is  one-one with dense
range, the resolvent  $R(\lambda,A)$ is defined  and bounded
for     $\lambda=re^{i\theta}$      where     $r>0$      and
$\omega<|\theta|\le   \pi$   and   satisfies   an   estimate
$\|\lambda      R(\lambda,A)\|\le      C_{\sigma}$       for
$\omega<\sigma\le |\theta|.$

Suppose  $(A_1,\ldots,A_m)$   is  a   commuting  family   of
sectorial operators where $A_k$ is of type $\omega_k$ for $1\le
k\le   m$,   and   let  $\omega=(\omega_1,\ldots,\omega_m)$.
Define   the   resolvent   for   $|\arg\lambda|>\omega$   by
$R(\lambda,A_1,\ldots,A_m)=\prod_{k=1}^mR(\lambda_k,A_k).$
Let ${\mathcal  A}$ be  the closed  subalgebra of ${\mathcal
L}(X)$  of  all  operators  $T$  so  that  $T$ commutes with
$R(\lambda,A_k)$  for  every  $k$  and  every $\lambda$ with
$|\arg\lambda|>\omega_k.$

If $F\in {\mathcal H}(\omega,{\mathcal  A})$ is of the  form
$F(z) = \prod_{k=1}^m(\lambda_k-z_k)^{-p_k}S$ where  $p_k\in
\mathbb   N\cup\{0\}$   and   $S\in\mathcal   A$  we  define
$F(A_1,\ldots,A_m)  =\prod_{k=1}^m   R(\lambda_k,A_k)^{p_k}S$
and then this definition can be extended by linearity to the
linear  span  of  such  functions,  which  we  call the {\it
rational  functions}, denoted  $\mathcal  R(\omega,\mathcal  A),$  in
$\mathcal H(\omega,{\mathcal A}).$

To  extend  this  definition  further  we  use the following
device.  Consider the algebra of all $(F,F(A_1,\ldots,A_m))$
for $F\in  \mathcal R(\omega,{\mathcal  A})$ as  a subset of
$\mathcal H(\omega,{\mathcal  A})\times \mathcal  A.$ Denote
by $\tau^*$ the sequential convergence $(F_n,T_n)\to  (F,T)$
if $F_n\to F (\tau)$  and $T_n\to T$ in  the strong-operator
topology.   Let ${\mathcal  B}$ be  the $\tau^*$-closure  of
this  set  (i.e.  the  smallest  set  which  is closed under
sequential convergence and contains  it).  Notice that  this
construction might involve  taking infinitely many iterations of  sequential
limits,  but  our  construction  actually  shows  that   two
iterations suffice.  It is  clear that ${\mathcal B}$ is  an
algebra.   Our next  task is  to show  that if $F\in\mathcal
H(\omega,\mathcal  A)$  there  is  at  most  one  choice  of
$T\in\mathcal  A$  so  that  $(F,T)\in\mathcal B$, this will
enable us to define $F(A_1,\ldots,A_m)$ unambiguously.

Consider      the      function      on      $\mathbb     C$
\begin{equation}\label{varphidef}              \varphi_n(z)=
\frac{n}{n+z}-\frac1{1+nz}\end{equation} and then define  on
$\mathbb C^m,$  $\psi_n(z)= \prod_{k=1}^m\varphi_n(z_k)$  so
that  $\psi_n  \in  H^{\infty}_0(\Sigma_{\sigma})$ for every
$\sigma>0.$  Then  $$\psi_n(A_1,\ldots,A_m)=   \prod_{k=1}^m
(\frac1n    R(-\frac1n,A_k)-n    R(-n,A_k))=V_n$$    is   an
approximate identity in the sense that $\sup \|V_n\|<\infty$
and $V_nx\to x$ for every $x\in X.$

If  $F\in  \mathcal  H(\omega,{\mathcal  A})$  then if $F\in
H^{\infty}  (\Sigma_{\sigma},{\mathcal  A})$  we  can define
\begin{equation}\label{int0} L_{n}(F)x= \left(\frac{-1}{2\pi i}
\right)^m      \int_{\Gamma_{\nu}}       \psi_n(\zeta)F(\zeta)
R(\zeta,A_1,\ldots,A_m)x\, d\zeta,\end{equation} as long  as
$\omega<\nu<\sigma.$ (Note that we are  using
short-hand and this is really a multiple contour  integral.)
An application  of Cauchy's  Theorem shows  that $L_{n}$  is
independent of the choice of  $\nu.$ By the Lebesgue  Dominated
Convergence  Theorem  $L_{n}:    \mathcal H(\omega,{\mathcal
A})\to {\mathcal A}$ is $\tau-$continuous if ${\mathcal  A}$
is equipped with the strong-operator topology.

If  $F$  is  rational  then  we  have  by a standard contour
integration,    \begin{equation}\label    {int}   L_{n}(F)x=
F(A_1,\ldots,A_m)V_nx\qquad  x\in  X.\end{equation}  Now the
map $(F,T)\to L_{n}(F)-T$ is continuous for $\tau^*$ and the
strong-operator   topology.       We   conclude   that    if
$(F,T)\in\mathcal B,$ $$  L_{_n}(F)x= TV_nx\qquad x\in  X.$$
Since $V_nx\to x$ for all  $x\in X,$ this shows that  $T$ is
uniquely determined by $F.$  Hence we can define  ${\mathcal
H}(A_1,\ldots,A_m;\mathcal   A)$   to   be   the   set    of
$F\in\mathcal H(\omega,{\mathcal A})$ such that for some $T$
we  have  $(F,T)\in{\mathcal  B}$  and  then  we  can define
$T=F(A_1,\ldots,A_m)$        for        $F\in       \mathcal
H(A_1,\ldots,A_m,\mathcal   A).$    The   space    $\mathcal
H(A_1,\ldots,A_m;\mathcal  A)$  is  an  algebra  and   $F\to
F(A_1,\ldots,A_m)$ is  an algebra  homomorphism.   For $F\in
\mathcal          H(A_1,\ldots,A_m;\mathcal           A)\cap
H^{\infty}(\Sigma_{\sigma};
\mathcal A)$  and  $\sigma>\nu>\omega$  then
(\ref{int0})   and   (\ref{int})   can   be   rewritten  as:
\begin{equation}\label{int2}     F(A_1,\ldots,A_m)V_nx     =
\left(\frac{-1}{2\pi i}\right)^m \int_{\Gamma_\nu}\psi_n(\zeta)
F(\zeta) R(\zeta,A_1,\ldots,A_m)x\,d\zeta.\end{equation}

If $F\in H^{\infty}_0(\Sigma_{\sigma};\mathcal A)$ then  the
integrals in (\ref{int2}) converge as $n\to\infty$.  We  can
show  by  approximating  the  integral  by Riemann sums that
$F\in\mathcal  H(\omega,\mathcal  A)$  and  then  we   have:
\begin{equation}\label{int3}      F(A_1,\ldots,A_m)x       =
\left(\frac{-1}{2\pi  i}\right)^m  \int_{\Gamma_\nu}   F(\zeta)
R(\zeta,A_1,\ldots,A_m)x\,d\zeta         \qquad         x\in
X.\end{equation}

It now follows that if $F\in\mathcal H(\omega,{\mathcal A})$
then  $(\psi_kF)\in\mathcal  H(A_1,\ldots,A_m;\mathcal  A)$ for each $k\in\mathbb N.$
Furthermore     if     $F_n\to     F(\tau)$     we      have
$(\psi_kF_n)(A_1,\ldots,A_m)\to  (\psi_kF)(A_1,\ldots,A_m)$   in
the strong-operator topology for each fixed $k.$.  From this
it follows that if $F_n\in\mathcal H(A_1,\ldots,A_m;\mathcal
A)$     and     $\sup\|F_n(A_1,\ldots,A_m)\|<\infty$     then
$F\in\mathcal H(A_1,\ldots,A_m)$ and $F_n(A_1,\ldots,A_m)\to
F(A_1,\ldots,A_m)$ in  the strong-operator  topology (indeed
we  have  convergence  on  each  $V_nx$).   In particular it
follows that $F\in\mathcal H(A_1,\ldots,A_m;\mathcal A)$  if
and only if $\sup_n\|(\psi_nF)(A_1,\ldots,A_m)\|<\infty.$

If   we   consider   the   scalar   functions  in  $\mathcal
H(A_1,\ldots,A_m) \subset\mathcal  H(A_1,\ldots,A_m;\mathcal
A)$ then we have defined the {\it joint functional  calculus}
for $(A_1,\ldots,A_m).$ We recall that a single operator $A$
has  {\it   an  $H^{\infty}(\Sigma_{\sigma})-$calculus}   if
$H^{\infty}(\Sigma_{\sigma})\subset   \mathcal   H(A).$  The
collection    $(A_1,\ldots,A_m)$    has    {\it    a   joint
$H^{\infty}(\Sigma_{\sigma})-$calculus         }          if
$H^{\infty}(\Sigma_{\sigma})\subset                 \mathcal
H(A_1,\ldots,A_m).$

\section{Rademacher-boundedness    and    related     ideas}
\label{Rademacher} \setcounter{equation}{0}

We recall (\cite{CPSW},\cite{We}) that a family $\mathcal F$
of bounded operators  on a Banach  space $X$ is  called {\it
Rademacher-bounded or R-bounded} with R-boundedness constant
$C$ if letting  $(\epsilon_k)_{k=1}^{\infty}$ be a  sequence
of independent  Rademachers on  some probability  space then
for  every  $x_1,\ldots,x_n  \in  X$  and $T_1,\ldots,T_n\in
\mathcal F$ we  have:  \begin{equation}\label{rad}  (\mathbb
E\|\sum_{k=1}^n\epsilon_k   T_kx_k\|^2)^{\frac12}   \le    C
(\mathbb    E\|\sum_{k=1}^n\epsilon_k    x_k\|^2)^{\frac12}.
\end{equation} It is important to note that this  definition
and the associated constant $C$ are unchanged if we  require
$T_1, \ldots,T_n$  to be  distinct in  (\ref{rad}) (see e.g.
\cite{CPSW}, Lemma 3.3).  The same remark applies to each  of
the following definitions.

We will introduce  two related weaker  notions.  Let  us say
that  $\mathcal  F$  is  {\it  weakly  Rademacher-bounded or
WR-bounded} with  WR-boundedness constant  $C$ if  for every
$x_1,\ldots,x_n  \in  X,  \  x_1^*,\ldots, x^*_n\in X^*$ and
$T_1,\ldots,T_n\in     \mathcal      L(X)$     we      have:
\begin{equation}\label{wrad}  \sum_{k=1}^n  |\langle T_kx_k,
x_k^*\rangle  |  \le  C  (\mathbb E(\|\sum_{k=1}^n\epsilon_k
x_k\|^2)^{\frac12}    (\mathbb    E(\|\sum_{k=1}^n\epsilon_k
x_k^*\|^2)^{\frac12}.  \end{equation}

Finally we  say that  $\mathcal F$  is {\it  unconditionally
bounded or U-bounded} with U-boundedness constant $C$ if for
every $x_1,\ldots,x_n\in X,\ x_1^*,\ldots,x_n^*\in X^*$  and
$T_1,\ldots,    T_n    \in     \mathcal    F$    we     have
\begin{equation}\label{uncb}            \sum_{k=1}^n|\langle
T_kx_k,x_k^*\rangle|\le      C      \max_{\epsilon_k=\pm1}\|
\sum_{k=1}^n\epsilon_k         x_k\|\max_{\epsilon_k=\pm1}\|
\sum_{k=1}^n\epsilon_k x^*_k\|.\end{equation}

The following Lemma is recorded for future reference:

\begin{Lem}\label{radetc} Let  $\mathcal F$  be a  subset of
$\mathcal  L(X)$.     Then   for  $\mathcal   F,$  R-bounded
$\Rightarrow$ WR-bounded  $\Rightarrow$ U-bounded.   If  $X$
has nontrivial Rademacher type then WR-bounded $\Rightarrow$
R-bounded.  \end{Lem}

We note that the only  really non-trivial part of the  Lemma
is the last sentence  and this follows easily  from Pisier's
characterization of spaces with non-trivial type as those in
which the Rademacher projection is bounded \cite{Pi2}.

We  shall  also  need  some  related  Banach space concepts.
Suppose          $(\epsilon_k)_{k=1}^{\infty}$           and
$(\eta_k)_{k=1}^{\infty}$  are   two  mutually   independent
sequences  of  Rademachers.    We  say that $X$ has property
$(\alpha)$ (see \cite{Pi1} and  \cite{LLLeM}) if there is  a
constant $C$ so that  for any $(x_{jk})_{j,k=1}^n\subset X$  and
any    $(\alpha_{jk})_{j,k=1}^n\subset\mathbb    C$    we   have
\begin{equation}\label{alpha}                       (\mathbb
E\|\sum_{j=1}^n\sum_{k=1}^n   \alpha_{jk}   \epsilon_j\eta_k
x_{jk}\|^2)^{\frac12} \le C\max_{j,k}|\alpha_{jk}|  (\mathbb
E\|\sum_{j=1}^n\sum_{k=1}^n                 \epsilon_j\eta_k
x_{jk}\|^2)^{\frac12}.  \end{equation}  We say that  $X$ has
property (A) \cite{LLeM} if there is a constant $C$ such that
for   any    $(x_{jk})_{j,k=1}^n\subset   X$    and   for    any
$(x^*_{jk})_{j,k=1}^n\subset        X^*$        we         have:
\begin{equation}\label{A} \sum_{j=1}^n\sum_{k=1}^n  |\langle
x_{jk},x_{jk}^*\rangle|        \le        C         (\mathbb
E\|\sum_{j=1}^n\sum_{k=1}^n                 \epsilon_j\eta_k
x_{jk}\|^2)^{\frac12}  (\mathbb  E\|\sum_{j=1}^n\sum_{k=1}^n
\epsilon_j\eta_k      x^*_{jk}\|^2)^{\frac12}.\end{equation}
Clearly $(\alpha)$ implies (A) and the converse holds if $X$
has  nontrivial  Rademacher  type;  this  is a fairly simple
deduction from the  boundedness of the  Rademacher projection
\cite{Pi2}.    Any  subspace   of  a  Banach  lattice   with
nontrivial cotype has  property $(\alpha)$ while  any Banach
lattice  has  property  (A).     It  is  also   observed  in
\cite{LLLeM}  that  $L_1/H_1$  has  $(\alpha).$ The Schatten
ideals $\mathcal C_p$ when  $1\le p\le \infty$ fail  to have
(A).

We shall say that $X$ has property $(\Delta)$ if there is  a
constant  $C$  so  that  for  any  $(x_{jk})_{j,k=1}^n\in X$
\begin{equation}\label{Delta}                       (\mathbb
E\|\sum_{j=1}^n\sum_{k=1}^j      \epsilon_j\eta_k
x_{jk}\|^2)^{\frac12}         \le         C         (\mathbb
E\|\sum_{j=1}^n\sum_{k=1}^n                 \epsilon_j\eta_k
x_{jk}\|^2)^{\frac12}.    \end{equation}  It  is  clear that
$(\Delta)$ is a  weaker property than  $(\alpha).$ It is  in
fact shared by all spaces with (UMD) and even analytic  UMD.
We recall  (\cite{Bl}) that  $X$ has  {\it analytic  UMD} if
every  $L_1-$bounded  analytic  martingale has unconditional
martingale differences.

\begin{Prop}\label{AUMD} Suppose $X$ has analytic UMD.  Then
$X$ has property ($\Delta$).\end{Prop}

\begin{proof}  Let  $(\tilde\epsilon_k)_{k=1}^{\infty}$  and
$(\tilde\eta_k)_{k=1}^{\infty}$ be two mutually  independent
sequences   of   Steinhaus    variables   (i.e.   each    is
complex-valued  and  uniformly   distributed  on  the   unit
circle).      By   applying   the  unconditionality  of  the
Rademachers  and  the  Khintchine-Kahane  inequality  it  is
sufficient to show the existence  of a constant $C$ so  that
for   any   $(x_{jk})_{j,k=1}^n$   we   have:    $$  \mathbb
E\|\sum_{j=1}^n\sum_{k=1}^j     \tilde\epsilon_j\tilde\eta_k
x_{jk}\|   \le    C   \mathbb    E\|\sum_{j=1}^n\sum_{k=1}^n
\tilde\epsilon_j\tilde\eta_k  x_{jk}\|.$$  To  see  this  we
define $f_j$ for $1\le  j\le 2n-1$ by $f_{2r-1}=  \sum_{j\le
r}\sum_{k\le   r}\tilde\epsilon_j\tilde\eta_k   x_{jk}$  and
$f_{2r}=              \sum_{j\le              r+1}\sum_{k\le
r}\tilde\epsilon_j\tilde\eta_k  x_{jk}.$  Let  $f_0=0.$ Then
$(f_j)$  is  an  analytic  martingale  and so for a suitable
constant $C$  depending only  on   $X$   we   have:     $$  \mathbb
E\|\sum_{r=0}^{n-1}(f_{2r+1}-f_{2r})\|     \le      C\mathbb
E\|f_{2n-1}\|.$$      This      yields      the      desired
inequality.\end{proof}

Since any  space with  (UMD) has  analytic (UMD)  this shows
that (UMD)-spaces have  ($\Delta$); actually a  direct proof
using  Rademacher  in  place  of  Steinhaus variables in the
above argument is possible for this case.  Thus the Schatten
classes $\mathcal C_p$ have  property $(\Delta)$ as long  as
$1<p<\infty.$  However  Haagerup  and  Pisier \cite{HP} show
that $\mathcal C_1$ (which has cotype 2) fails analytic  UMD
and  their  argument   actually  shows  it   fails  property
$(\Delta)$.      This   implies   that   $C(K)$-spaces    of
infinite dimension also fail ($\Delta$) since $\mathcal C_1$
is finitely representable in any such space.

We  now  come  to  an  important  result  relating the above
properties to Rademacher-boundedness.  Some similar  results
are shown in \cite{CPSW}.

\begin{Thm}\label{uncseries} Suppose  $(U_k)_{k=1}^{\infty}$
and $(V_k)_{k=1}^{\infty}$ are two sequences of operators in
${\mathcal      L}(X)$       satisfying      $$       \sup_n
\sup_{\epsilon_k=\pm1} \|  \sum_{k=1}^n \epsilon_k  U_k\|\le
M<\infty$$   and   $$   \sup_n   \sup_{\epsilon_k=\pm1}   \|
\sum_{k=1}^n \epsilon_k V_k\|\le M<\infty.$$ Suppose further
$\mathcal F\subset \mathcal L(X)$  is a family of  operators
which  commutes  with  each  $U_k$  and  each  $V_k$  and is
R-bounded with constant $R.$ Then:  \begin{enumerate}  \item
The  sequence   $(U_n)_{n=1}^{\infty}$  is   R-bounded  with
constant  $M.$  \item  If  $X$  has  property $(\alpha)$ the
collection $\{\sum_{k=1}^n\alpha_k T_kU_kV_k:  \ n\in\mathbb
N,\ |\alpha_k|\le 1,\ T_k\in\mathcal F\} $ is R-bounded with
constant  $CRM^2$  where  $C$   depends  only  on  $X.$   In
particular  the  family  $\{\sum_{k=1}^n  \alpha_k  U_kV_k:\
n\in\mathbb  N,\  |\alpha_1|,\ldots,|\alpha_n|\le  1\}$   is
R-bounded with  constant $CM^2$  where $C$  depends only  on
$X.$  \item  If  $X$  has  property  $(A)$  then  the family
$\{\sum_{k=1}^n    \alpha_k    U_kV_k:\    n\in\mathbb   N,\
|\alpha_1|,\ldots,|\alpha_n|\le  1\}$  is  WR-bounded   with
constant $CM^2$ where $C$ depends only on $X.$ \item If  $X$
has property $(\Delta)$ then the set  $\{\sum_{k=1}^nU_kV_k:
\ n\in\mathbb N\}$ is  R-bounded with constant $CM^2$  where
$C$ depends only on $X.$ \end{enumerate} \end{Thm}

\begin{proof} (1)  We use  the remark  that it  is enough to
establish     (\ref{rad})     for     distinct     operators
$T_1,\ldots,T_n.$ If $x_1,\ldots,x_n\in X$ and $\alpha_k=\pm
1$          then          $$\mathbb          E((\sum_{k=1}^n
\epsilon_kU_k)(\sum_{k=1}^n               \epsilon_k\alpha_k
x_k))=\sum_{k=1}^n   \alpha_k    U_kx_k$$   and    hence   $$
\sup_{\alpha_k=\pm1}\|\sum_{k=1}^{\infty}\alpha_kU_kx_k\|
\le    M     (\mathbb    E(\|     \sum_{k=1}^n    \epsilon_k
x_k\|^2)^{\frac12}.$$ This  proves (1)  and indeed  a rather
stronger result.

(2)  Let  $S_j  =\sum_{k=1}^{\infty}\alpha_{jk}T_{jk}U_kV_k$
where  $T_{jk}\in\mathcal  F$   and  $(\alpha_{jk})$  is   a
finitely  nonzero   collection  of   complex  numbers   with
$|\alpha_{jk}|\le  1.$  Suppose  $x_1,\ldots,x_n\in  X.$  We
first note that \begin{align*} \|\sum_{k=1}^nU_kV_kx_k\|  &=
\|     \mathbb     E\left(     (\sum_{k=1}^n      \epsilon_k
U_k)(\sum_{k=1}^n  \epsilon_k  V_k  x_k)\right)\|  \\ &\le M
(\mathbb        E\|\sum_{k=1}^n        \epsilon_k        V_k
x_k)\|^2)^{\frac12}.\end{align*} We will  also use the  fact
(Lemma 3.13  of \cite{CPSW})  that there  is a  constant $C$
depending only on $X$ so that for $(x_{jk})_{j,k=1}^n\in  X$
we   have    from   property    $(\alpha)$,   $$    (\mathbb
E_{\epsilon}\mathbb       E_{\eta}\|\sum_{j=1}^n\sum_{k=1}^n
\alpha_{jk}   \epsilon_j\eta_kT_{jk}x_{jk}\|^2)^{\frac12}\le
CR               (\mathbb                E_{\epsilon}\mathbb
E_{\eta}\|\sum_{j=1}^n\sum_{k=1}^n
\epsilon_j\eta_kx_{jk}\|^2)^{\frac12}.$$

Hence                \begin{align*}                 (\mathbb
E_{\epsilon}\|\sum_{j=1}^n\epsilon_jS_jx_j\|^2)^{\frac12} &=
(\mathbb
E_{\epsilon}\|\sum_{k=1}^{\infty}U_kV_k\sum_{j=1}^n\alpha_{jk}
\epsilon_jT_{jk}x_j\|^2)^{\frac12}   \\   &\le   M  (\mathbb
E_{\epsilon}\mathbb
E_{\eta}\|\sum_{j=1}^n\sum_{k=1}^{\infty}\alpha_{jk}
\epsilon_j\eta_k   V_kT_{jk}x_j\|^2)^{\frac12}\\   &\le  CRM
(\mathbb                                 E_{\epsilon}\mathbb
E_{\eta}\|\sum_{j=1}^n\sum_{k=1}^{\infty}   \epsilon_j\eta_k
V_kx_j\|^2)^{\frac12}\\       &\le       CRM^2      (\mathbb
E_{\epsilon}\|\sum_{j=1}^n                        \epsilon_j
x_j\|^2)^{\frac12}.\end{align*} This proves (2).

(3)  Let  $S_j=\sum_{k=1}^{\infty}\alpha_{jk}U_kV_k$   where
$(\alpha_{jk})$   is   a   finitely   nonzero   matrix  with
$|\alpha_{jk}|\le 1.$ In this case if $x_1,\ldots,x_n\in  X$
and    $x_1^*,\ldots,x_n^*\in    X^*$    we    note    that:
\begin{align*}  \sum_{j=1}^n  |\langle  S_jx_j,x_j^*\rangle|
&\le        \sum_{j=1}^n\sum_{k=1}^{\infty}         |\langle
V_kx_j,U_k^*x^*_j\rangle|     \\     &\le     C     (\mathbb
E_{\epsilon}\mathbb               E_{\eta}                \|
\sum_{j=1}^n\sum_{k=1}^{\infty}             \epsilon_j\eta_k
V_kx_j\|^2)^{\frac12} (\mathbb E_{\epsilon}\mathbb  E_{\eta}
\|     \sum_{j=1}^n\sum_{k=1}^{\infty}      \epsilon_j\eta_k
U_k^*x^*_j\|^2)^{\frac12}\\   &\le   CM^2   (\mathbb   E  \|
\sum_{j=1}^n  \epsilon_jx_j\|^2)^{\frac12}  (\mathbb  E   \|
\sum_{j=1}^n \epsilon_jx^*_j\|^2)^{\frac12} \end{align*}

(4) We use  the proof of  (2).  This  time we again  use the
fact  it   suffices  to   consider  the   operators  without
repetition.   So we  consider $S_j=  \sum_{k=1}^jU_kV_k$ and
repeat the proof of (2) with $\alpha_{jk}=1$ if $k\le j$ and
$0$ otherwise  and replace  each $T_{jk}$  by the  identity.
Using  (\ref{Delta})  in  place  of  (\ref{alpha}) gives the
desired conclusion.  \end{proof}

We conclude this  section with a  useful Lemma.   In fact in
the  case  of  R-boundedness,   this  result  is  found   in
\cite{We}.

\begin{Lem} \label{powers2} Suppose $0<\sigma<\pi$ and $F\in
H^{\infty}(\Sigma_{\sigma},   {\mathcal   L}(X)).$   Suppose
$0<\sigma_0<\nu<\sigma$ and for some $M<\infty$ and $a>1$, and for each
$t\in\mathbb      R$      the      set       $\{F(a^kte^{\pm
i\nu})\}_{k\in\mathbb   Z}$   is   U-bounded  (respectively,
WR-bounded; respectively, R-bounded) with constant bounded by
$M$  (independent  of   $t).$  Then  the   family  $\{F(z):\
z\in\Sigma_{\sigma_0}\}$   is   U-bounded,    (respectively,
WR-bounded; respectively, R-bounded).\end{Lem}

\begin{proof} We give the  proof in the U-boundedness  case,
the others  being similar.   We  first make  the observation
that    it    suffices    to    consider   the   case   when
$\nu=\frac{\pi}{2}$  as  one  can  make  the  transformation
$z=w^{2\nu/\pi}.$ In this case  we have the formula  $$ F(z)
=\frac1\pi  \int_{-\infty}^{\infty}  F(it)  \Re  (z-it)^{-1}
dt.$$  We  write  $$  F_1(z)=\frac{1}{\pi} \int_{0}^{\infty}
F(it)  \Re  (z-it)^{-1}  dt$$  and  $$  F_2(z)=\frac{1}{\pi}
\int_{0}^{\infty}  F(-it)  \Re  (z+it)^{-1}  dt  $$  so that
$F(z)=F_1(z)+F_2(z).$

Note that for  a suitable constant  $C$ we have  an estimate
$0\le \Re (z\pm  it)^{-1} \le C|z||\min  (t^{-2}, |z|^{-2})$
whenever $z\in \Sigma_{\sigma_0}.$

Now  suppose  $x_1,\ldots,x_n\in  X,\  x_1^*,\ldots,x_n^*\in
X^*.$ Suppose $z_1,\ldots, z_n\in \Sigma_{\sigma_0}.$ Let us
suppose  that  $m_j\in  \mathbb   Z$  are  chosen  so   that
$a^{m_j}\le |z_j|\le  a^{m_j+1}.$ We  have:   \begin{align*}
\sum_{j=1}^n    |\langle    F_1(z_j)x_j,x_j^*\rangle|   &\le
\frac1{\pi}       \sum_{j=1}^n       \int_0^{\infty}|\langle
F(ia^{m_j}t)x_j,x_j^*\rangle|   \Re   (z_j-ia^{m_j}t)^{-1}
a^{m_j}dt     \\     &\le     \frac{aC}{\pi}    \sum_{j=1}^n
\int_0^{\infty}|\langle F(ia^{m_j}t)x_j,x_j^*\rangle| \min
(1,t^{-2})          dt          \\          &\le          C'
\max_{\epsilon_j=\pm1}\|\sum_{j=1}^n\epsilon_jx_j\|
\max_{\epsilon_j=\pm1}\|\sum_{j=1}^n\epsilon_jx_j^*\|\end{align*}
for a suitable constant $C'.$ A similar argument can be done
for $F_2.$\end{proof}

\section{Functional calculus  for operator-valued  functions
}\label{operator}\setcounter{equation}{0}

Let  us  suppose  $A$  is  sectorial  of  type  $\omega$ and
$\sigma>\omega.$ We  let $\mathcal  A$ denote  as in Section
\ref{functional} the algebra of all bounded operators  which
commute with $A.$

Before we prove our basic estimate for an operator-valued functional calculus, we will describe in Lemma \ref{unc} and Proposition \ref{Hinfty} the connection between the $H^{\infty}-$calculus and unconditional expansions in the underlying Banach space.

\begin{Lem}\label{unc}   Suppose   that   $A$   admits    an
$H^{\infty}(\Sigma_{\sigma})$-calculus,   and   that   $f\in
H^{\infty }_0(\Sigma_{\sigma}).$ Then there is a constant $C$ so that
for  any  $t>0$  and   any  finitely  nonzero  sequence   $(
\alpha_k)_{k\in  \mathbb  Z}$  we  have:    $$  \|\sum_{k\in
\mathbb   Z}\alpha_kf(2^ktA)\|   \le   C\max_{k\in   \mathbb
Z}|\alpha_k|.$$ Furthermore for every $x\in X$ and $t>0$ the
series  $   \sum_{k\in\mathbb  Z}   f(2^ktA)x  $   converges
unconditionally in $X.$ \end{Lem}

\begin{proof}    We    can    assume   $\max_{k\in   \mathbb
Z}|\alpha_k|\le  1.$  For  a  suitable  constants $C,C'$ and
$\epsilon>0$  we  have  \begin{align*}  \|\sum_{k\in \mathbb
Z}\alpha_kf(2^ktA)\|&\le    C\sup_{z\in\Sigma}    \sum_{k\in
\mathbb Z}|f(2^kz|\\  &\le CC'\sup_{z\in  \Sigma} \sum_{k\in
\mathbb
Z}\left(\frac{2^k|z|}{1+2^{2k}|z|^2}\right)^{\epsilon}\end{align*}
and the last quantity is finite.

For the last part observe that for any bounded sequence $(\alpha_k)_{k\in\mathbb  Z}$  and
$t>0,$ the  series $\sum_{k\in
Z}\alpha_k  f(2^ktA)  x$  must  converge  to  $g(A)x$  where
$g(z)=\sum_{k\in           Z}\alpha_k            f(2^ktz)\in
H^{\infty}(\Sigma_{\sigma}).$ \end{proof}

\begin{Prop}\label{intrep}           Suppose           $F\in
H_0^{\infty}(\Sigma_{\sigma},{\mathcal  A}).$  Then  for any
$\omega<\nu<\sigma$,    $0<s<1,$    and    any   $x\in   X,$
\begin{equation}\label{altdef}      F(A)x=       \frac{-1}{2\pi
i}\int_{\Gamma_{\nu}}\zeta^{-s}F(\zeta)A^sR(\zeta,A)x\,
d\zeta.  \end{equation} \end{Prop}

\begin{proof} First note that $A^sR(\lambda,A)$ is a bounded
operator for $\lambda\in \Gamma_{\nu}$ which is given by the
integral    $$     A^s    R(\lambda,A)x     =    \frac{-1}{2\pi
i}\int_{\Gamma_{\nu'}}\zeta^{s}(\lambda-\zeta)^{-1}R(\zeta,A)x\,
d\zeta,  $$  if  $\omega<\nu'<\nu.$  This  gives an estimate
$\|A^sR(\lambda,A)\|\le C_s |\lambda|^{s-1}$ and shows  that
the  integral  in  (\ref{altdef})  converges  as  a  Bochner
integral.    It  is  clear  that  we only need establish the
formula  if  $x=\varphi_n(A)y$  (see  (\ref{varphidef})) for
some    $y\in     X.$    To     do    this     we    compute
\begin{align*}F(A)\varphi_n(A)x   &=   F(A)\varphi_n^2(A)y\\
&=(A^s\varphi_n(A))              (F(A)A^{-s}\varphi_n(A))y\\
&=\frac1{2\pi                            i}(A^s\varphi_n(A))
\int_{\Gamma_{\nu}}\zeta^{-s}\varphi_n(\zeta)F(\zeta)
R(\zeta,A)y      \,d\zeta\\      &=      \frac1{2\pi      i}
\int_{\Gamma_{\nu}}\zeta^{-s}\varphi_n(\zeta)F(\zeta)(A^s
\varphi_n(A))R(\zeta,A)y\,d\zeta\\    &=    \frac1{2\pi   i}
\int_{\Gamma_{\nu}}\zeta^{-s}\varphi_n(\zeta)F(\zeta)(A^s
R(\zeta,A))x\,d\zeta\\ \end{align*} Now using the  Dominated
Convergence Theorem we obtain (\ref{altdef}).\end{proof}

Let us rewrite (\ref{altdef}) by using the  parameterization
$\zeta=|t|e^{i(\text{sgn }t)\nu}$  for $-\infty<t<\infty.$  We
introduce         the         notation        $h^{\rho}_s(z)
=z^{s}(e^{i\rho}-z)^{-1}.$       Then       for        $F\in
H^{\infty}_0(\Sigma_{\sigma},\mathcal       A)$        where
$\sigma>\nu>\omega,$   \begin{align*}   F(A)x&=  \frac1{2\pi
i}\int_{-\infty}^{\infty}e^{i(1-s)(\text{sgn
}t)\nu}|t|^{-s}F(|t|e^{i(\text{sgn                 }t)\nu})A^s
R(|t|e^{i(\text{sgn             }t)\nu},A)x\,             dt\\
&=\frac{e^{i(1-s)\nu}}{2\pi        i}        \int_0^{\infty}
F(te^{i\nu})h_s^{\nu}(t^{-1}A)x\,\frac{dt}{t}     \\&\     +
\frac{e^{-i(1-s)\nu}}{2\pi        i}         \int_0^{\infty}
F(te^{-i\nu})h_s^{-\nu}(t^{-1}A)x\,\frac{dt}{t}\end{align*}
This      can       then      be       reformulated      as:
\begin{equation}\label{basic1}      F(A)x=       \frac1{2\pi
i}\int_1^2\left(            M_+(t)+M_-(t)\right)\frac{dt}{t}
\end{equation}     where      \begin{equation}\label{basic2}
M_{\pm}(t)=     e^{\pm     i(1-s)\nu}      \sum_{k\in\mathbb
Z}F(2^{-k}t^{-1}e^{\pm           i\nu})h_s^{\pm\nu}(2^ktA)x.
\end{equation}

We   first   make    an   essentially   trivial    deduction
characterizing the $H^{\infty}$-calculus.

\begin{Prop}\label{Hinfty}Suppose $\nu>\omega$ and  $0<s<1$.
Consider  the   conditions:     \begin{equation}\label{unc1}
\sup_{t>0}\sup_N   \sup_{\epsilon_k=\pm1}    \|\sum_{k=-N}^N
\epsilon_k    (2^{k}t)^{(1-s)}A^s    R(2^kte^{\pm   i\nu},A)
\|<\infty,\end{equation} Then (\ref{unc1}) is necessary  for
$A$ to  admit an  $H^{\infty}(\Sigma_{\sigma})-$calculus for
some  $\sigma<\nu$  and  sufficient  for  $A$  to  admit  an
$H^{\infty}(\Sigma_{\sigma})$-calculus       for       every
$\sigma>\nu.$ \end{Prop}

\begin{proof}  Necessity  follows  immediately  from   Lemma
\ref{unc} for the  functions $h_s^{\pm\nu}.$ Conversely,  by
(\ref{unc1}),  if  $f\in  H^{\infty}(\Sigma_{\sigma})$ where
$\sigma>\nu$ we obtain by (\ref{basic1}) and  (\ref{basic2})
$$  \|(\varphi_nf)(A)\|  \le  C$$  independent  of $n.$ This
implies that $f\in{\mathcal H}(A).$ \end{proof}

Our  main  result  is  also  easy  from  (\ref{basic1})  and
(\ref{basic2}).

\begin{Thm}\label{multiplier}Suppose    $A$    admits     an
$H^{\infty}(\Sigma_{\sigma})-$calculus       and       $F\in
H^{\infty}(\Sigma_{\rho};{\mathcal     A})$     for     some
$\rho>\sigma.$  Suppose  further  that  the set $\{F(z):z\in
\Sigma_{\rho}\}$  is  $U$-bounded.    Then  $F\in  {\mathcal
H}(A,{\mathcal A}).$\end{Thm}

{\it Remarks.} (1) Of course the Theorem holds if we  assume
the stronger property that $\{F(z):z\in \Sigma_{\rho}\}$  is
WR-bounded or R-bounded.

(2) For  Hilbert spaces  and certain  operators on  $L_2(X)$
such an operator-valued  functional calculus is  constructed
in  \cite{LLLeM}  Theorem  5.2  and  \cite{LLeM}.  These are
cases   when   the   U-boundedness   condition  is  satisfied
automatically.  See also \cite{FM}.  In \cite{CP} there are
constructions based on  transference results which  work for
generators of  bounded $c_0$-groups  on UMD-spaces  and some
other special cases.

\begin{proof} As before we consider $\varphi_nF=F_n$ so that
$F_n\in{\mathcal H}(A,{\mathcal  A}).$ It  suffices to  show
$\sup\|F_n(A)\|<\infty.$  Referring  to  (\ref{basic1})  and
(\ref{basic2}) with some fixed $0<s<1$ and $\rho>\nu>\sigma$
for  $x\in  X,\  x^*\in  X^*$  with $\|x\|,\|x^*\|\le 1$, we
obtain  the  estimate  for  $1\le  t\le  2$:    $$  |\langle
M_{\pm}(t)x,x^*\rangle|\le     \sum_{k\in\mathbb     Z}|\langle
F_n(2^{-k}t^{-1}e^{\pm                      i\nu})g(2^ktA)x,
g(2^ktA)^*x^*\rangle|      $$      where     $g(z)=(h_s^{\pm
\nu}(z))^{\frac12}.$  Suppose  $C$  is  the  $U$-boundedness
constant of $\{F(z):z\in\Sigma_{\sigma}\}.$ Then $$ |\langle
M_{\pm}(t)x,x^*\rangle|                \le                 C
\sup_{\epsilon_k=\pm1}\sup_N\|\sum_{|k|\le
N}\epsilon_kg(2^ktA)\|^2.$$ Hence by Lemma \ref{unc} we have
$$\sup_n\|F_n(A)\| <\infty.$$ \end{proof}

Let us apply this to the case of two commuting operators:

\begin{Thm}\label   {Two}   Suppose   $A,B$   are  commuting
sectorial    operators,    such    that    $A$    admits   a
$H^{\infty}(\Sigma_{\sigma})-$calculus                   and
$\omega(B)<\sigma'.$              Suppose              $f\in
H^{\infty}(\Sigma_{\rho}\times\Sigma_{\sigma'})$       where
$\sigma<\rho<\pi$    is    such    that   $\{f(w,\cdot):w\in
\Sigma_{\rho}\}$ is contained in ${\mathcal H}(B).$  Suppose
further   the   set   $\{f(w,B):w\in   \Sigma_{\rho}\}$   is
U-bounded.  Then $f\in {\mathcal H}(A,B)$ (i.e.  $f(A,B)$ is
a bounded operator).  \end{Thm}

\begin{proof} We  define $F(w)=f(w,B)$  and note  that $F\in
H^{\infty}(\Sigma_{\rho};\mathcal A);$  this follows  easily
from   the   integral   representation   (\ref{int}).    Our
conditions   and   Theorem   \ref{multiplier}   ensure  that
$F\in{\mathcal H}(A;\mathcal  A).$ It  is only  necessary to
check  that  this  implies  $f\in{\mathcal  H}(A,B)$  and of
course  $F(A)=f(A,B).$  But   this  follows  directly   from
(\ref{int}),      (\ref{int2})      and      the     remarks
thereafter.\end{proof}

{\it Example.} Let us show by example that Theorem \ref{Two}
is  close  to  the  best  possible.   Let $B$ be a sectorial
operator on  $X$.   Suppose $0<\sigma<\pi$  and consider the
space       $L_2(\{-1,1\}^{\Sigma_{\sigma}};X)$        where
$\{-1,1\}^{\Sigma_{\sigma}}$ has the usual product  measure.
Denote  by  $\epsilon_z$  the  co-ordinate  maps  for  $z\in
\Sigma_{\sigma}$.  Let Rad $X$ denote the closed linear span
of   the   functions    $\{\epsilon_z\otimes   x:       z\in
\Sigma_{\sigma},\ x\in X\}.$ We define $\tilde B=I\otimes B$
on $L_2(X)$ and restrict it to the subspace Rad $X$ which is
invariant.    We  define  $A$  on Rad $X$  by  $$  A(\sum_{z\in
\Sigma_{\sigma}}\epsilon_zx_z)=   \sum_{z\in\Sigma_{\sigma}}
z\epsilon_zx_z$$    with    domain    consisting    of   all
$\sum\epsilon_z x_z\in L_2$ so that $\sum  z\epsilon_zx_z\in
L_2.$

Clearly  $A$  has  an $H^{\infty}(\Sigma_{\sigma})-$calculus
and  $f(A,\tilde B)$  is  bounded   if  and  only  if   the  family
$\{f(z,B):z\in\Sigma_{\sigma}\}$ is R-bounded.

We remark that the  reader who prefers separable  spaces can
easily modify this example when $X$ is separable to replace Rad $X$ by a separable subspace (just take a dense countable subset of $\Sigma_{\sigma}.$

\section{R-boundedness   and   the   functional  calculus  }
\label{props} \setcounter{equation}{0}

We now consider strengthenings of the boundedness conditions
in the definition of sectoriality.   Let $A$ be a  sectorial
operator  and  let  $\omega(A)$  denote  the  infimum of all
$\sigma$ so that $A$ is  of type $\sigma.$ We will  say that
$A$   is   {\it   R-sectorial,  (respectively  WR-sectorial,
respectively U-sectorial)} if there exists $0<\sigma<\pi$ so
that     the     family     of     operators      $\{\lambda
R(\lambda,A):|\arg\lambda|>\sigma\}$      is       R-bounded
(respectively WR-bounded, respectively U-bounded).  We  then
define   $\omega_R(A)$,   (respectively    $\omega_{WR}(A)$,
respectively $\omega_U(A)$)  to be  the infimum  of all such
$\sigma.$  We  will  say  $A$ is {\it $H^{\infty}-$sectorial
(respectively,     $RH^{\infty}-$sectorial,     respectively
$WRH^{\infty}-$sectorial)} if there exists a  $0<\sigma<\pi$
so   that   $A$   admits   an  $H^{\infty}(\Sigma)-$calculus
(respectively,        such        that        the        set
$\{f(A):\|f\|_{H^{\infty}(\Sigma_{\sigma})}\le    1\}$    is
R-bounded,   respectively   such   that   the  set  $\{f(A):
\|f\|_{H^{\infty}(\Sigma_{\sigma})}\le 1\}$ is  WR-bounded).
The infimum  of all  such $\sigma$  is denoted $\omega_H(A)$
(respectively         $\omega_{RH}(A)$,         respectively
$\omega_{WRH}(A)$).

There are certain obvious and trivial relationships  between
these concepts.   Clearly  R-sectorial implies  WR-sectorial
implies U-sectorial and whenever these concepts are defined,
$\omega_R(A)\ge \omega_{WR}(A)\ge \omega_U(A)\ge \omega(A).$
Similarly       $RH^{\infty}-$       sectorial       implies
$WRH^{\infty}-$sectorial implies $H^{\infty}-$sectorial  and
$\omega_{RH}(A)\ge     \omega_{WRH}(A)\ge     \omega_H(A)\ge
\omega(A).$

We now turn to less trivial observations:

\begin{Prop}\label      {UH}       Suppose      $A$       is
$H^{\infty}-$sectorial    and    $U$-sectorial.         Then
$\omega_H(A)\le \omega_U(A).$\end{Prop}

\begin{proof}    Let     us    assume     that    $\{\lambda
R(\lambda,A):|\arg  \lambda|\ge  \nu\}$  is  U-bounded  with
constant $K$ where  $\nu>\omega(A)$, and that  $\sigma>\nu.$
We      will      show      that      $A$      admits     an
$H^{\infty}(\Sigma_{\sigma})-$calculus.  We use  Proposition
\ref{Hinfty}.  Fix  some $0<s<1.$ We  can assume that  there
exists        $\rho>\sigma$        so        that         $$
\sup_N\sup_{\epsilon_k=\pm1}\sup_{t>0}\|\sum_{|k|\le
N}\epsilon_kh_s^{\pm\rho}(A)\| =M<\infty$$  and so  that $A$
admits  an  $H^{\infty}(\Sigma_{\tau})-$calculus  for   some
$\tau<\rho.$

Now suppose $x\in X$ and $x^*\in X^*.$ Then for any $N$  and
$\epsilon_j=\pm1$   we   have   $$   |\langle   \sum_{|k|\le
N}\epsilon_k h_s^{\nu}(2^ktA)x,x^*\rangle |\le M\|x\|\|x^*\|
+            \sum_{|k|\le            N}             |\langle
(h_s^{\nu}(2^ktA)-h_s^{\rho}(2^ktA))x,x^*\rangle|.$$ By  the
resolvent                                          equation,
$$h_s^{\nu}(2^ktA)-h_s^{\rho}(2^ktA)=(e^{i(\rho-\nu)}-1)2^{-k}t^{-1}e^{i\nu}
R(2^{-k}t^{-1}e^{i\nu},A) h_s^{\rho}(A).$$ Since $A$ has  an
$H^{\infty}(\Sigma_{\tau})-$calculus    we    can     define
$g(z)=(h_s^{\rho}(z))^{\frac12}$    and    note    that   $$
\sup_N\sup_{\epsilon_k=\pm1}\|\sum_{|k|\le      N}\epsilon_k
g(2^ktA)\|\le C$$ where $C$ is independent of $t.$ Thus,  by
the    U-boundedness    of    $\{\lambda   R(\lambda,   A):\
\arg\lambda=\nu\}$    $$     \sum_{|k|\le    N}     |\langle
2^{-k}t^{-1}R(2^{-k}t^{-1}e^{i\nu},A)g(2^ktA)x,g(2^ktA)^*x^*\rangle|
\le  KC   \|x\|\|x^*\|.$$  It   follows  that   $$  |\langle
\sum_{|k|\le N}\epsilon_k h_s^{\nu}(2^ktA)x,x^*\rangle  |\le
(M+2KC)\|x\|\|x^*\|$$      and      this      gives       $$
\sup_N\sup_{\epsilon_k=\pm1}\sup_{t>0}\|\sum_{|k|\le
N}\epsilon_kh_s^{\nu}(A)\| \le M+2KC<\infty.$$ Combined with
a similar estimate for $-\nu$ we obtain the result by  using
Proposition \ref{Hinfty}.  \end{proof}

In  order  to  study  an  analytic  semigroup with generator
$(-A)$  it   is  of   particular  interest   to  know   that
$\omega_H(A)<\frac{\pi}{2}.$  Therefore  we  use Proposition
\ref{UH} to improve on a result in \cite{HP2}.

\begin{Cor}\label{improve} Let  $(-A)$ generate  an analytic
contractive and positive semigroup on $L_p(\Omega,\mu)$  for
some                   $1<p<\infty.$                    Then
$\omega_H(A)<\frac{\pi}{2}.$\end{Cor}

\begin{proof}    It    is    shown    in   \cite{HP2}   that
$\omega_H(A)<\pi$  and  in   \cite{We2},  Section  5,   that
$\omega_R(A)<\frac{\pi}{2}.$ Hence we can apply  Proposition
\ref{UH}.\end{proof}

We remark  that it  is an  open problem  (cf.   \cite{CDMY})
whether    $\omega_H(A)=\omega(A)$    whenever    $A$     is
$H^{\infty}-$sectorial.  The next Theorem gives some results
in this direction.

\begin{Thm}\label{RadHinfty}    Suppose     $A$    is     an
$H^{\infty}-$sectorial  operator  on  a  Banach  space  $X$.
Then:     \begin{enumerate}  \item   If  $X$   has  property
$(\alpha)$   then   $A$   is   $RH^{\infty}-$sectorial   and
$\omega_H(A)=\omega_{RH}(A)=\omega_R(A)=\omega_U(A).$  \item
If    $X$    has     property    $(A)$    then     $A$    is
$WRH^{\infty}-$sectorial                                 and
$\omega_H(A)=\omega_{WRH}(A)=\omega_{WR}(A)=\omega_U(A).$
\item If $X$ has property $(\Delta)$ then $A$ is R-sectorial
and  $\omega_H(A)=\omega_R(A)=\omega_U(A).$  \end{enumerate}
\end{Thm}

\begin{proof}    (1)    Assume    that    $A$    admits   an
$H^{\infty}(\Sigma_{\sigma})$-calculus.              Suppose
$\sigma<\nu<\pi.$  Suppose  $0<s<1$  and  let   $g_{\pm}(z)=
(h_s^{\pm \nu}(z))^{\frac12}.$  We then  can argue  by Lemma
\ref{unc}      that      $$     \sup_N\sup_{\epsilon_k=\pm1}
\|\sum_{k=-N}^N \epsilon_k  g_{\pm}(2^ktA)\| \le  M<\infty$$
independent of $t.$ Hence by Lemma \ref{uncseries}   the
family $\{\sum_{|k|\le N}\alpha_k h_s^{\pm \nu}(2^ktA)\}$ is
R-bounded with constant bounded  independent of $t.$ Now  by
(\ref{basic1})  and  (\ref{basic2})   it  follows  that   if
$\sigma'>\nu$                                           then
$\{f(A):\|f\|_{H^{\infty}(\Sigma_{\sigma'})}\le   1\}$    is
Rademacher-bounded.            Indeed      for       $f_k\in
H^{\infty}_0(\Sigma_{\sigma'})$  and  $x_k\in  X$  for $1\le
k\le     n,$     we     have     \begin{align*}    &(\mathbb
E\|\sum_{k=1}^n\epsilon_kf_k(A)x_k\|^2)^{\frac12}   \le   \\
&\le 4  \max_\pm \sup_{t>0}\sup_{N\in\mathbb  N} (\mathbb  E
\|\sum_{k=1}^n    \epsilon_k(\sum_{|j|\le    N}   f_n(e^{\pm
i\nu}t^{-1}2^{-j})h_s^{\pm
i\nu}(2^jtA))x_k\|^2)^{\frac12}.\end{align*} It follows that
$\omega_{RH}(A)=\omega_H(A).$  Now  clearly  $\omega_U(A)\le
\omega_R(A)\le  \omega_{RH}(A)$  and  so  (1)  follows  from
Proposition \ref{UH}.

(2) is very similar and we omit it.

(3) Here we use Lemma \ref{powers2}.  Suppose $A$ admits  an
$H^{\infty}(\Sigma_{\sigma})-$calculus      and      suppose
$\sigma'>\nu>\sigma.$  We  show  that  the  sequence $\{2^kt
R(2^kte^{\pm i\nu}):\ k\in\mathbb Z\}$ is Rademacher-bounded
with constant independent of $t.$ To do this we note that if
$N_1>N_2$              $$              2^{N_1}tR(2^{N_1}te^{
i\nu},A)-2^{N_2}tR(2^{N_2}te^{i\nu},A)=-\sum_{j=N_2+1}^{N_1}t2^{j-1}AR(2^jte^{
i\nu},A)     R(2^{j-1}te^{i\nu},A).$$     Let     $k(z)    =
z(e^{i\nu}-z)^{-1}(e^{i\nu}-2z)^{-1}.$                   Let
$u(z)=(k(z))^{\frac12}\in  H^{\infty}(\Sigma_{\sigma}).$  We
observe                       that                        $$
\sup_{N_1>N_2}\sup_{\epsilon_j=\pm1}\|\sum_{j=N_2+1}^{N_1}
\epsilon_j u(2^{-j}t^{-1}A)\| \le M<\infty$$ independent  of
$t$  by  Lemma  \ref{unc}.    Applying Lemma \ref{uncseries}
yields  that  $$  \{\sum_{j=N_2+1}^{N_1}  k(2^{-j}t^{-1}A):\
N_1>N_2\}$$ is Rademacher-bounded with constant  independent
of       $t.$       But       this       implies        that
$\{2^{N_1}tR(2^{N_1}te^{i\nu},A)-2^{N_2}tR(2^{N_2}te^{i\nu},A):
\  N_1>N_2\}$  is  also  Rademacher-bounded  with   constant
independent  of  $t$  and   hence  (taking  limits)  so   is
$\{2^ntR(2^nte^{i\nu},A):\  n\in\mathbb   Z\}.$  A   similar
argument   for   $-\nu$   and   an   application   of  Lemma
\ref{powers2}   shows   that   $\omega_R(A)\le  \nu.$  Hence
$\omega_R(A)\le \omega_H(A).$  The proof  is finished  as in
(1).\end{proof}

As a Corollary  to the proof  of Theorem \ref{RadHinfty}  we
obtain some  additional information  on the  operator-valued
calculus considered in Theorem \ref{multiplier}:

\begin{Cor}\label{4.4}   Assume   that   $X$   has  property
$(\alpha)$ and let $\mathcal  F\subset \mathcal L(X)$ be  an
R-bounded set.   If  $A$ is  $H^{\infty}$-sectorial then for
any    $\sigma>\omega_H(A)$    the    set   $\{F(A):\   F\in
H^{\infty}(\Sigma_{\sigma},\mathcal A),\ F(\zeta)\in\mathcal
F\ \forall \zeta\in\Sigma_{\sigma}\}$ is R-bounded.\end{Cor}

\begin{proof} Adapt the proof of Theorem \ref{RadHinfty} (1)
using    the    fact    that    the    set  $$\{\sum_{|k|\le
N}T_kh_s^{\pm\nu}(2^ktA):    T_k\in\mathcal  F\cap  \mathcal
A\}$$      is       R-bounded,      again       by      Lemma
\ref{uncseries}.\end{proof}

\section{The joint $H^{\infty}$-calculus and sums of  closed
operators}\setcounter{equation}{0}

First we consider the joint functional calculus.

\begin{Thm} \label{joint} Suppose $A$ and $B$ are  commuting
$H^{\infty}-$sectorial   operators   such   that   $B$    is
$WRH^{\infty}-$sectorial.            Then      for       any
$\sigma>\omega_{H}(A)$  and  $\sigma'>\omega_{WRH}(B)$   the
pair $(A,B)$  has a  joint $H^{\infty}(\Sigma_{\sigma}\times
\Sigma_{\sigma'})-$calculus.  \end{Thm}

\begin{proof}   We   need   only   observe   that  if  $f\in
H^{\infty}(\Sigma_{\sigma}\times \Sigma_{\sigma'})$ then the
family  $\{f(z,B):    z\in H^{\infty}(\Sigma_{\sigma})\}$ is
WR-bounded        and        then        apply       Theorem
\ref{multiplier}.\end{proof}

We can now apply Theorem \ref{RadHinfty} to obtain a  result
of  Lancien,  Lancien  and  Le  Merdy \cite{LLLeM} (see also
\cite{AFM}).  Note that their argument depends on the  quite
technical discretization developed in \cite{FM}.

\begin{Cor} \label{joint2}  (Lancien, Lancien  and Le  Merdy
\cite{LLLeM}) If $X$  has property (A)  then if $A$  and $B$
are  commuting  $H^{\infty}-$sectorial  operators,  for  any
$\sigma>\omega_{H}(A)$ and $\sigma'>\omega_{H}(B)$ the  pair
$(A,B)$   has   a   joint  $H^{\infty}(\Sigma_{\sigma}\times
\Sigma_{\sigma'})-$calculus.  \end{Cor}

If $A$ and $B$ are commuting sectorial operators on a Banach
space $X$ with  $\omega(A)+\omega(B)<\pi$ then one  can show
that the  closure $\overline{A+B}$  of $A+B$  with $\mathcal
D(A+B)=\mathcal  D(A)\cap\mathcal   D(B)$  is   a  sectorial
operator  with  $\omega(A+B)\le   \max(\omega(A),\omega(B))$
(see \cite{DA}).  However, in many applications one needs to
show   that   $A+B$   is   already   closed   on   $\mathcal
D(A)\cap\mathcal D(B).$ We now give a criterion for this.

\begin{Thm}\label{sumthm} Suppose $A$ and $B$ are  commuting
sectorial operators such that $A$ is  $H^{\infty}-$sectorial
and  $B$  is  R-sectorial and $\omega_H(A)+\omega_R(B)<\pi.$
Then   $A+B$   is    closed   on   the    domain   $\mathcal
D(A)\cap\mathcal D(B),$  there is a constant $C$ such that
\begin{equation}\label{sum}   \|Ax\|+\|Bx\|\le    C\|Ax+Bx\|
\qquad x\in\mathcal D(A)\cap\mathcal D(B)\end{equation}  and
$(A+B)$ is invertible  if either $A$  or $B$ is  invertible.
Furthermore if  $X$ has  property $(\alpha)$  then $A+B$  is
again        R-sectorial        with       $\omega_R(A+B)\le
\max(\omega(A),\omega(B)).$\end{Thm}

{\it  Remarks.}  (1)  Let  us  compare this theorem with the
well-known Dore-Venni  Theorem.   It is  shown in  \cite{DV}
that (\ref{sum}) holds if $X$ is a UMD-space and $A,B$  both
have bounded imaginary  powers (BIP), with  $$ \|A^{is}\|\le
Ce^{\theta_A|s|},\quad   \|B^{is}\|\le    Ce^{\theta_B|s|}$$
where    $\theta_A+\theta_B<\pi.$     In    a     UMD-space,
R-sectoriality  is  weaker  than  (BIP)  (see \cite{CP}) and
$H^{\infty}-$sectoriality is  stronger than  (BIP).   But in
many  applications  $A$  is  an  operator  known  to have an
$H^{\infty}$-calculus,  e.g.    $A=-\Delta$  or  $A=d/dt$  on
$L_p(X)$ where $X$ is UMD and $1<p<\infty.$ Thus the  weaker
assumption on  $B$ does  lead to  more general  results, see
e.g.  Theorem \ref{maxreg} below.

(2) Some special cases  of Theorem \ref{sumthm} where  shown
in \cite{We}  Theorem 5.2,  \cite{We2} and  more recently in
\cite{CP}  (e.g.  if  $A$  is  the  generator  of a strongly
continuous group in a UMD-space).

(3) An  extension to  non-commuting sums  will be  given in a
forthcoming paper \cite{Str}.

\begin{proof}     Choose     $\sigma,     \sigma'$      with
$\omega_H(A)<\sigma,$       $\omega_R(B)<\sigma'$        and
$\sigma+\sigma'<\pi.$ The  function $f(w,z)=w(w+z)^{-1}$  is
in  $H^{\infty}(\Sigma_{\sigma}\times\Sigma_{\sigma'})$  and
the  set  $f(w,B)=-wR(-w,B)$ for $w\in\Sigma_{\sigma}$  is  an  R-bounded  family.
  Applying  Theorem  \ref{Two} we have
$f\in\mathcal H(A,B).$ We  can see this  implies (\ref{sum})
either by applying Proposition 2.7 in \cite{LLeM} or by  the
following simple direct  argument based on  our construction
of  the  functional  calculus.    Defining $\varphi_n$ as in
(\ref{varphidef})   we   note   that  $A\varphi_n(A)^2$  and
$B\varphi_n(B)^2$     are     bounded     operators    since
$z\varphi_n(z)^2\in  H^{\infty}_0(\Sigma_{\tau})$  for   any
$\tau<\pi.$ Now if  $x\in\mathcal D(A)\cap\mathcal D(B)$  we
have      $$       f(A,B)(A+B)\varphi_n(A)^2\varphi_n(B)^2x=
\varphi_n(A)^2\varphi_n(B)^2Ax.$$          Thus           $$
\|\varphi_n(A)^2\varphi_n(B)^2Ax\|\le
C\|\varphi_n(A)^2\varphi_n(B)^2(A+B)x\|$$              where
$C=\|f(A,B)\|.$ Letting $n\to\infty$ yields the result.

Now   assume   $X$    has   property   $(\alpha)$.       For
$\max(\sigma,\sigma')<\rho<\pi$  and  $|\arg\mu|\ge  \rho  $
consider  the  functions  $f_{\mu}(w,z)=\mu(\mu-w-z)^{-1}\in
H^{\infty}(\Sigma_{\sigma}\times\Sigma_{\sigma'}).$     Note
that                    $$                     f_{\mu}(w,B)=
\frac{\mu}{\mu-w}((\mu-w)R(\mu-w,B)).$$                Since
$\mu(\mu-w)^{-1}$  is  bounded  uniformly  for $|\arg\mu|\ge
\rho$ and  $w\in\Sigma_{\sigma}$ and  also $|\arg(\mu-w)|\ge
\sigma'$ this  collection of  operators in  R-bounded.  Now,
Corollary \ref{4.4}  yields that  the set  $f_{\mu}(A,B)=\mu
R(\mu,A+B)$ is R-bounded for $|\arg\mu|\ge \rho.$\end{proof}

Applying Theorem \ref{RadHinfty} gives:

\begin{Cor}  \label  {sumthm2}  Suppose  $X$  has   property
$(\Delta)$ (e.g. if  $X$ has analytic  (UMD)).  Suppose  $A$
and $B$ are commuting $H^{\infty}-$sectorial operators  such
that $\omega_H(A)+\omega_H(B)<\pi.$ Then $A+B$ is closed  on
the domain $\mathcal D(A)\cap\mathcal D(B)$ and  (\ref{sum})
holds.  \end{Cor}

{\it Example.} We now show by example that both  Corollaries
\ref{joint2} and  \ref{sumthm2} are  nearly optimal.   To do
this we let $(\epsilon_j)$  and $(\eta_k)$ be as  before two
sequence  of  mutually   independent  Rademachers  on   some
probability   space   $(\Omega,\mathbb   P).$   We    define
$\text{Rad}_2(X)$  to  be  the  subspace  of $L_2(\Omega;X)$
spanned  by  functions  of  the form $\epsilon_j\eta_kx$ for
$j,k\in\mathbb N$  and $x\in  X.$ Let  $A$ be  defined by $$
A(\sum_{j,k}\epsilon_j\eta_kx_{jk})=\sum_{j,k}(2j+1)!
\epsilon_j\eta_kx_{jk}$$ with the natural domain and let  $$
B(\sum_{j,k}\epsilon_j\eta_kx_{jk})=\sum_{j,k}(2k)!
\epsilon_j\eta_kx_{jk}$$ with its natural domain.

  Both  $A$  and  $B$  are  $H^{\infty}-$sectorial   with
$\omega_H(A)=\omega_H(B)=0.$ Clearly $f\in{\mathcal H}(A,B)$
if  and   only  if   $$  \sum_{j,k}\epsilon_j\eta_kx_{jk}\to
\sum_{j,k} f((2j+1)!,(2k)!)  \epsilon_j\eta_k x_{jk}$$ is  a
bounded operator.

If          $(A,B)$           has          a           joint
$H^{\infty}(\Sigma_{\sigma}\times\Sigma_{\sigma'})-$calculus
then  since  the  map  $f\to  (f((2j+1)!,(2k)!))_{j,k}$ maps
$H^{\infty}(\Sigma_{\sigma}\times  \Sigma_{\sigma'})  $ onto
$\ell_{\infty}(\mathbb  N^2)$  we  must  have  that  $X$ has
property $(\alpha)$ (and hence so does $\text{Rad}_2(X).$)

If $f(w,z)=w(w+z)^{-1}$ defines a bounded operator then by a
limiting         argument         the         map         $$
\sum_{j,k}\epsilon_j\eta_kx_{jk}       \to        \sum_{k\le
j}\epsilon_j\eta_kx_{jk}$$ is also  bounded, i.e.   $X$ (and
$\text{Rad}_2(X)$) has property $(\Delta).$\qed

Let  us  note  that  we  can  recapture  the  main result of
\cite{We}  on  maximal  regularity.    Suppose  $-A$  is the
generator     of     an     analytic     semigroup      i.e.
$\omega(A)<\frac{\pi}2.$   Then   $A$   has  maximal  $L_p-$
regularity   for   $1<p<\infty$   if   the   Cauchy  problem
\begin{equation}\label{cauchy}       y'(t)+Ay(t)=f(t),\qquad
t>0,\quad   y(0)=0\end{equation}   has   for   every   $f\in
L_p(\mathbb R_+,X)$  a solution  $y:\mathbb R_+\to  X$ which
satisfies  the  estimate  \begin{equation}\label{regularity}
\|y'\|_{L_p(\mathbb  R_+,X)}+\|Ay\|_{L_p(\mathbb  R_+,X)}\le
C\|f\|_{L_p(\mathbb R_+,X)}.\end{equation}  If we  denote by
$\tilde B$  the derivative  $d/dt$ on  $\tilde X=L_p(\mathbb
R_+,X)$  and  by  $\tilde  A$  the extended operator $\tilde
Af(t)=Af(t)$ then (\ref{regularity})  is equivalent to:   $$
\|\tilde  Ay\|_{\tilde   X}  +\|\tilde   By\|_{\tilde  X}\le
C\|(\tilde A+\tilde  B)y\|_{\tilde X}.$$  Thus we  can apply
Theorem \ref{sumthm}:

\begin{Thm}\label{maxreg} Suppose $X$ is a Banach space with
(UMD) and suppose $A$ is an R-sectorial operator on $X$ with
$\omega_R(A)<\frac{\pi}2.$    Then    $A$    has     maximal
$L_p$-regularity for $1<p<\infty.$\end{Thm}

\begin{proof} Since $X$ has  (UMD), we have that  $\tilde B$
is      $H^{\infty}-$sectorial      and     $\omega_H(\tilde
B)=\frac{\pi}{2}$ (see e.g.  \cite{HP2}).  It is easy to see
that  $\omega_R(\tilde  A)=\omega_R(A)$  so  the result is a
consequence of Theorem \ref{sumthm}.\end{proof}

{\it   Remarks.}   (1)   If   $X$   has  (UMD)  and  $A$  is
$H^{\infty}-$sectorial    then     $\omega_H(A)<\frac{\pi}2$
implies maximal $L_p$-regularity by Corollary \ref{sumthm2}.

(2) It is shown in \cite{KL} that any non-Hilbertian  Banach
space  with  an  unconditional  basis  admits  a   sectorial
operator  $A$   with  $\omega(A)<\frac{\pi}2$   but  failing
maximal $L_p-$regularity.

(3)   It   is   shown   in   \cite{We}  that  the  condition
$\omega_R(A)<\frac{\pi}2$       actually       characterizes
$L_p$-maximal regularity.  The operator approach to  maximal
regularity    has     a    long     history,    see     e.g.
\cite{DA},\cite{DO}, \cite {G} and \cite{We2}.

\section{$L_1$-spaces   and   $C(K)$-spaces}   \label{Groth}
\setcounter{equation}{0}

We recall that a GT-space is a Banach space $X$ so that  the
Grothendieck  theorem  is  valid,  i.e. if $T:X\to\ell_2$ is
bounded  then  $T$  is  absolutely  summing and for some $C$
independent of $T$, $\pi_1(T)\le C\|T\|$; see \cite{Pi0} for
a   full   discussion.      Examples   of  such  spaces  are
$L_1-$spaces,  their  quotients  by  reflexive  spaces   and
$L_1/H_1$ ( \cite{Bo1}, \cite{BD} and \cite{Pi0}).

\begin{Prop}\label{GT} Suppose $X$ is a GT-space of cotype 2
and that  $A$ is  a $H^{\infty}-$sectorial  operator on $X.$
Then if  $\omega_H(A)<\nu$ and  $0<s<1$ there  is a constant
$C$ so  that if  $x\in X,$  \begin{equation}\label{absolute}
C^{-1}\|x\|\le   \int_{\Gamma_{\nu}}   \|   A^sR(\zeta,A)x\|
\frac{|d\zeta|}{|\zeta|^s}    \le    C\|x\|    \qquad   x\in
X.\end{equation} \end{Prop}

\begin{proof} This is just  a version of a  classical result
of Lindenstrauss and Pe\l czy\'nski \cite{LP} on  uniqueness
of     unconditional     bases     in     $\ell_1.$      Let
$g(z)=(h_s^{\nu}(z))^{\frac12}.$  Then  for  some   constant
$C_0$ independent of $x,t$ we have if $|\alpha_k|\le 1$  for
$|k|\le N,$  $$ \|  \sum_{|k|\le N}\alpha_k  g(2^ktA)x\| \le
C_0\|x\| \qquad x\in X, t>0$$ by Lemma \ref{unc}. Now  since
$X$ has cotype 2, there exists a constant $C_1$  independent
of    $x,t,$    so    that    \begin{equation}\label{cotype}
\left(\sum_{k\in{\mathbb                                 Z}}
\|g(2^ktA)x\|^2\right)^{\frac12} \le  C_1 \|x\|  \qquad x\in
X, t>0.\end{equation}  Now suppose  $x\in X$  and $t>0$  and
choose  by  the  Hahn-Banach  theorem  $x_k^*\in  X^*$  with
$\|x_k^*\|=1$                  and                  $\langle
h_s^{\nu}(2^ktA)x,x_k^*\rangle=\|h_s^{\nu}(2^ktA)x\|.$
Consider  the  map  $S_x:X\to\ell_2(\mathbb  Z)$  defined by
$S_xy=(\langle g(2^ktA)y,x_k^*\rangle)_{k\in\mathbb Z}.$  By
(\ref{cotype}) $S$  is bounded  and $\|S_x\|\le  C_1.$ Hence
for some constant $C_2$ we have $\pi_1(S_x)\le C_2.$

Now \begin{align*}  \sum_{|k|\le N}\|h_s^{\nu}(2^ktA)x\|  &=
\sum_{|k|\le  N}  \langle   h_s^{\nu}(2^jtA)x,x_k^*\rangle\\
&\le   \sum_{|k|\le   N}\|   S_x   g(2^ktA)x\|\\   &\le  C_2
\sup_{|\alpha_k|\le       1}\|\sum_{|k|\le        N}\alpha_k
g(2^ktA)x\|\\  &\le  C_2  C_0  \|x\|.    \end{align*}  If we
integrate  for  $1\le  t\le  2$  we obtain the right-half of
(\ref{absolute})  (cf.    (\ref{basic1}) and (\ref{basic2})).
The   left-half   follows    from   the   equation:       $$
x=\lim_{n\to\infty}\frac{-1}{2\pi
i}\int_{\Gamma_{\nu}}\varphi_n(\zeta)\zeta^{-s}A^s
R(\zeta,A)x\,d\zeta.$$ \end{proof}

The same argument yields:

\begin{Prop}\label{GT*}  Suppose  $X^*$  is  a  GT-space  of
cotype 2 and that  $A$ is a $H^{\infty}-$sectorial  operator
on $X.$  Then if  $\omega_H(A)<\nu$ and  $0<s<1$ there  is a
constant  $C$  so  that  if  $x\in  X,$  $$\frac1C   \|x^*\|\le
\int_{\Gamma_{\nu}}        \|         (A^sR(\zeta,A))^*x^*\|
\frac{|d\zeta|}{|\zeta|^s} \le C\|x^*\| \qquad x^*\in X^*.$$
\end{Prop}

{\it Remarks.}  Let us  point out  that Proposition \ref{GT}
implies  that  very  few  operators  on  $L_1$  can  have an
$H^{\infty}-$calculus.  This statement can be made much more
precise  but  since  the  techniques  required  are   rather
specialized we will defer this to a later paper and  instead
note  the  following  simple  application, which effectively
shows that no reasonable differential operator on $L_1$  can
have an $H^{\infty}-$calculus.

\begin{Prop} Suppose $X$ is a  GT-space of cotype 2 and  $A$
is an $H^{\infty}-$sectorial operator  on $X.$ If $Y$  is an
infinite-dimensional closed reflexive subspace of  $\mathcal
D(A)$ (with the graph norm) then $A$ is bounded on $Y$  (and
so $Y$ is closed in $X$).

\end{Prop}

\begin{proof} We use  the notation of  Proposition \ref{GT}.
In  particular  notice  that  (\ref{absolute})  implies  the
existence    of    an    isomorphic    embedding     $T:X\to
L_1(\Gamma_{\nu}, |d\zeta|;\  X)$ defined  by $$  Tx (\zeta)
=|\zeta|^{-\frac12}  A^{\frac12}R(\zeta,A)x.$$  Fix   $0\neq
\lambda\in\Gamma_{\nu}.$   Then   $R(\lambda,A)$   maps  $X$
isomorphicly onto $\mathcal  D(A)$ (with the  graph norm).
Let   $Y_0=   R(\lambda,   A)^{-1}Y$;   then   $Y_0$  is  an
infinite-dimensional reflexive subspace of $X$.  Then (since
$Y_0$ contains no copy of $\ell_1$) the set $\{\|Ty(\zeta)\|
; \ \|y\|\le 1,\ y\in Y_0\}$ is equi-integrable.
We show that this implies that  the
(bounded)  operator  $A^{\frac12}R(\lambda,  A)$ satisfies a
lower bound on $Y_0.$ Indeed if not there exists a  sequence
$(y_n)$    in    $Y_0$     so    that    $\|y_n\|=1$     and
$\|A^{\frac12}R(\lambda,A)y_n\|\to 0.$ But, by the resolvent
equation,   $\|A^{\frac12}R(\zeta,A)y_n\|\to   0$   for  all
$\zeta\in\Gamma_{\nu}\setminus\{0\}.$         Now         by
(\ref{absolute}) and equi-integrability, we have $\|y_n\|\to
0$  which  gives  a  contradiction.    Now applying the same
argument to  $Y_1=A^{\frac12}R(\lambda,A)Y_0$ gives  a lower
bound on $AR(\lambda,A)^2$ on $Y_0.$ Thus $R(\lambda,A)$ has
a  lower  bound  on  $Y_0.$  Since  $Y=R(\lambda,A)Y_0$ this
implies the result.  \end{proof}

Let us  note, with  respect to  differential operators, that
embeddings  of   Sobolev  spaces   into  $L_1(\Omega)$   are
Dunford-Pettis operators.

In view of this comparative rarity of $H^{\infty}-$sectorial
operators in this setting, it is not surprising that we  may
substantially improve  the results  of this  paper for these
special  spaces.    Our  first  deduction  is that for these
special spaces, Theorem \ref{Two} can be improved  by
removing the U-boundedness assumption.

\begin{Thm} Suppose $X$ is  a Banach space such  that either
$X$ or $X^*$  is a GT-space  of cotype 2.  Suppose $A,B$ are
commuting  sectorial  operators  such  that  $A$  admits  an
$H^{\infty}-$calculus.                  Suppose        $f\in
H^{\infty}(\Sigma_{\sigma}\times  \Sigma_{\sigma'})$   where
$\sigma>\omega_H(A)$ and  $\sigma'>\omega(B).$ Suppose  that
for each $z\in \Sigma_{\sigma}$, $f_z(w)=f(z,w)\in  \mathcal
H(B)$ and $$  \sup_{z\in\Sigma_{\sigma}}\|f(z,B)\|<\infty.$$
Then $f\in  \mathcal H(A,B)$  (i.e.   $f(A,B)$ is  bounded).
\end{Thm}

\begin{proof} Let us assume that $X$ is a GT-space of cotype
2, the other case  is similar.  By  Proposition \ref{intrep}
if  $f_n(w,z)=\varphi_n(w)\varphi_n(z)f(w,z)$  we  can write
\begin{align*}            f_n(A,B)x&=            \frac{-1}{2\pi
i}\int_{\Gamma_{\nu}}\zeta^{-s}f_n(\zeta,B)A^sR(\zeta,A)x\,d\zeta\\
&=                                               \frac{-1}{2\pi
i}\int_{\Gamma_{\nu}}\zeta^{-s}\varphi_n(\zeta)
f(\zeta,B)A^sR(\zeta,A)\varphi_n(B)x\,d\zeta    \end{align*}
and  this  leads  immediately to $\sup_n\|f_n(A,B)\|<\infty$
which implies the boundedness of $f(A,B).$\end{proof}

It is clear  now if $A$  is $H^{\infty}-$sectorial on  $L_1$
with $\omega_H(A)<\pi/2$ then one can apply the above result
to conclude  that $A$  has $L_1-$maximal  regularity.   More
generally  we  have  the  following  result  (suggested by a
question of Gilles Lancien):

\begin{Thm}\label{final} Suppose  either $X$  or $X^*$  is a
GT-space of  cotype 2  (e.g. if  $X=L_1,\ C(K)$  or the disk
algebra $A(\mathbb D)).$ If $A$ is an  $H^{\infty}-$sectorial
operator  on  $X$  with  $\omega_H(A)<\pi/2$  then  $A$  has
$L_p-$maximal regularity for $1<p<\infty.$\end{Thm}

\begin{proof} Let us prove this for $X$ a GT-space of cotype
2    as    the    other    case    is    dual.       Suppose
$\omega_H(A)<\nu<\pi/2$ and $0<s<1.$ If $u>0$ then for $x\in
X$ by Proposition \ref{intrep} $$ Ae^{-uA}x =  \frac{-1}{2\pi
i}\int_{\Gamma_{\nu}}\zeta^{1-s}e^{-u\zeta}A^s R(\zeta,A)x\,
d\zeta.$$  Now  suppose  $f\in  L_p(\mathbb  R,X).$  We will
estimate the norm  of $S_{\delta}$ where  $$ S_{\delta}f(v)=
\int_{\delta}^{\infty}Ae^{-uA} f(v-u)du.$$ If we let  $G(t)=
A^s R(|t|e^{i(\text{sgn }t)\nu},A)$  then we have  an estimate
$$\|Ae^{-uA}x\|                   \le                    C_0
\int_{-\infty}^{\infty}|t|^{1-s}e^{-cu|t|}\|G(t)x\|dt$$
where $c=\cos\nu>0$  and so  $$ \|S_{\delta}f(v)  \| \le C_0
\int_{-\infty}^{\infty}\int_{-\infty}^{v-\delta}|t|^{1-s}e^{c(u-v)|t|
}\|G(t)f(u) \|du\,dt.$$  Now if  $g\in L_q(\mathbb  R)$ with
$g\ge           0$           we           have            $$
\int_{u+\delta}^{\infty}e^{c(u-v)|t|}g(v)dv \le \frac1{c|t|}
\int_u^{\infty}(v-u)e^{c(u-v)t}\frac1{v-u}\int_u^vg(w)dw\,dv.$$
Hence    the    left-hand    side    is   estimated   by   $
C_1|t|^{-1}(\mathcal  Mg)(t)$  where  $\mathcal  M$  is  the
Hardy-Littlewood maximal function.

Hence $$ \int_{-\infty}^{\infty}g(v)\|S_{\delta}f(v)\|dv \le
C_0C_1   \int_{-\infty}^{\infty}\int_{-\infty}^{\infty}|t|^{-s}
\mathcal M g(u)\|G(t)f(u)\|  du\,dt .$$ However  Proposition
\ref{GT}   implies   that    we   have   an    estimate   $$
\int_{-\infty}^{\infty}|t|^{-s}\|G(t)f(u)\|dt            \le
C_2\|f(u)\|.$$     Substituting     in     we     have    $$
\int_{-\infty}^{\infty}g(v)\|S_{\delta}f(v)\|dv  \le  C_0C_1C_2
\int_{-\infty}^{\infty}{\mathcal M}g(u)\|f(u)\|du $$ and since $M$  is
bounded on  $L_q$ this  establishes a  uniform bound  on the
operators  $S_{\delta}.$  Letting  $\delta\to  0$ yields the
result.\end{proof}

\end{document}